\documentclass[aop,MSNbibl,number,citesort,dvips]{arximspdf}

\doi{10.1214/11-AOP712} 
\volume{41}
\issue{2}
\pubyear{2013}
\firstpage{1055}
\lastpage{1071}

\makeatletter
\newtheorem{theo}{Theorem}
\newtheorem{prop}{Proposition}
\newtheorem{coro}{Corollary}
\newtheorem{lema}{Lemma}
\def\ccA{\mathcal{A}}
\def\ccF{\mathcal{F}}
\def\ccM{\mathcal{M}}
\def\va{\varepsilon}
\def\ra{\rightarrow}
\def\E{\mathbf{E}}
\def\P{\mathbf{P}}
\def\R{\mathbb{R}}
\def\ls{\leq}
\def\gs{\geq}
\makeatother

\begin{document}
\begin{frontmatter}

\title{The complete characterization of a.s. convergence of orthogonal series}
\runtitle{Convergence of orthogonal series}

\begin{aug}
\author[A]{\fnms{Witold} \snm{Bednorz}\corref{}\thanksref{t1}\ead[label=e1]{wbednorz@mimuw.edu.pl}}
\thankstext{t1}{Supported in part by the Funds of Grant MENiN N N201 397437.}
\runauthor{W. Bednorz}
\affiliation{Warsaw University}
\address[A]{Department of Mathematics\\
Warsaw University\\
Banacha 2\\
Warsaw, 02-097\\
Poland\\
\printead{e1}} 
\end{aug}

\received{\smonth{10} \syear{2010}}
\revised{\smonth{8} \syear{2011}}

%
\begin{abstract}
In this paper we prove the complete characterization
of a.s. convergence of orthogonal series in terms
of existence of a majorizing measure. It means that for a given
$(a_n)^{\infty}_{n=1}$,
$a_n>0$, series $\sum^{\infty}_{n=1}a_n\varphi_n$ is a.e. convergent for
each orthonormal
sequence $(\varphi_n)^{\infty}_{n=1}$ if and only if there exists a
measure $m$ on
\[
T=\{0\}\cup\Biggl\{\sum^m_{n=1}a_n^2, m\gs1\Biggr\}
\]
such that
\[
\sup_{t\in T}\int^{\sqrt{D(T)}}_0 (m(B(t,r^2)))^{-{1}/{2}}\,dr<\infty,
\]
where $D(T)=\sup_{s,t\in T}|s-t|$ and $B(t,r)=\{s\in T\dvtx  |s-t|\ls r\}$.
The presented approach is based on weakly majorizing measures and a
certain partitioning scheme.
\end{abstract}

%

\end{frontmatter}

\section{Introduction}\label{sect0}

An orthonormal sequence $(\varphi_n)^{\infty}_{n=1}$ on a probability
space $(\Omega,\ccF,\P)$
is a sequence of random variables $\varphi_n\dvtx \Omega\ra\R$ such that
$\E\varphi_n^2=1$
and $\E\varphi_n\varphi_m=0$ whenever $n\neq m$. In this paper we
consider the question of how to characterize the sequences of
$(a_n)^{\infty}_{n=1}$ for which the series
\[
\sum^{\infty}_{n=1}a_n\varphi_n \mbox{ converges a.e. for any
orthonormal } (\varphi_n)^{\infty}_{n=1}
\]
on any probability spaces $(\Omega,\ccF,\P)$. Note that
we can assume $a_n>0$ for $n\gs1$. The answer is based on the analysis
of the set
\[
T=\Biggl\{ \sum^m_{n=1}a_n^2\dvtx  m\gs1\Biggr\}\cup\{0\}.
\]
The classical Rademacher--Menchov theorem (see~\cite{Kas,Mo-Ta})
states that\break
$\sum^{\infty}_{n=1}a_n^2 \times \log^2(n+1)$ suffices for $\sum^{\infty
}_{n=1}a_n\varphi_n$
convergence. Another well-known observation (see~\cite{Tan2}) is the
following theorem.
\begin{theo}\label{thm0}
For each orthonormal sequence $(\varphi_n)^{\infty}_{n=1}$ the series\break
$\sum^{\infty}_{n=1}a_n\varphi_n$ converges a.e. if and only if
\[
\E\sup_{m\gs1}\Biggl(\sum^m_{n=1}a_n\varphi_n\Biggr)^2<\infty.
\]
\end{theo}

The consequence of the above result is that the main problem can be
reformulated in terms of sample boundedness of all orthogonal processes
on $T$. We say that process $X(t)$, $t\in T$, is of orthogonal
increments if
%
\begin{equation}\label{ineq2}
\E\bigl(X(s)-X(t)\bigr)^2=|s-t|\qquad \mbox{for } s,t\in T.
\end{equation}
There is a bijection between orthonormal series $\sum^{\infty
}_{n=1}a_n\varphi_n$ and processes with orthogonal increments on $T$.
Namely for each sequence $(\varphi_n)^{\infty}_{n=1}$ we define processes
\[
X(t)=\sum^m_{n=1}a_n\varphi_n \qquad\mbox{for } t=\sum^m_{n=1}a_n^2, X(0)=0,
\]
and for each orthogonal process $X(t)$, $t\in T$, we define
the orthonormal sequence by
\[
\varphi_m=a_m^{-1}\Biggl(X\Biggl(\sum^m_{n=1}a_n^2\Biggr)-X\Biggl(\sum^{m-1}_{n=0}a_n^2\Biggr)\Biggr)
\qquad\mbox{for } m>1,
\]
and $\varphi_1=a_1^{-1}(X(a_1^2)-X(0))$. By Theorem~\ref{thm0}, each
orthogonal series\break $\sum^{\infty}_{n=1}a_n\varphi_n$
is a.e. convergent if and only if there exists a universal constant
$\ccM<\infty$ such that
%
\begin{equation}\label{square}
\E\sup_{t\in T}|X(t)-X(0)|^2\ls\ccM
\end{equation}
for all processes $X(t)$, $t\in T$ that satisfy (\ref{ineq2}).

We treat the generalized question and consider any
$T\subset\R$. The best tool which is used to study the sample boundedness
of orthogonal processes on $T$ are majorizing measures. Let
$B(t,r)=\{s\in T\dvtx  |s-t|\ls r\}$ and $D(T)=\sup_{s,t\in T}|s-t|$. We
say that a probability measure $m$ on $T$
is majorizing (in the orthogonal setting) if
\[
\sup_{t\in T}\int^{\sqrt{D(T)}}_0 (m(B(t,r^2)))^{-{1}/{2}}\,dr<\infty.
\]
We say that a process $X(t)$, $t\in T$, is of suborthogonal increments if
%
\begin{equation}\label{ineq3}
\E|X(s)-X(t)|^2\ls|s-t|\qquad \mbox{for } s,t\in T.
\end{equation}
Corollary 1 proved in~\cite{Bed3} states that
the existence of a majorizing measure is the necessary and sufficient
condition for the sample boundedness of all
suborthogonal processes.
Moreover by Theorem 3.2 in~\cite{Bed1} (see also~\cite{Tal1})
we have the following theorem.
\begin{theo}\label{thm2}
For each process $X(t)$, $t\in T$, that satisfies (\ref{ineq3}),
the following inequality holds:
\[
\E\sup_{s,t\in T}|X(s)-X(t)|^2\ls16\cdot5^{{5}/{2}}
\biggl(\sup_{t\in T}\int^{\sqrt{D(T)}}_0 (m(B(t,r^2)))^{-
{1}/{2}}\,dr\biggr)^2.
\]
\end{theo}

Consequently the existence of a majorizing measure is always sufficient
for the a.e. convergence of orthogonal series
$\sum^{\infty}_{n=1}a_n\varphi_n$.

The problem is that the class of orthogonal processes is
significantly smaller than the class of processes that verify (\ref
{ineq3}). Only recently
Paszkiewicz proved in~\cite{Pas1,Pas2}, using advanced methods of
entropy of interval, that the existence of a majorizing measure is also
necessary for all orthogonal processes to satisfy (\ref{square}). This
motivated our research for an alternative approach
entirely based on the generic chaining; see~\cite{Tal2,Tal3}.
We use the Fernique's idea of constructing a majorizing measure. We say
that a probability measure $\mu$ on $T$
is weakly majorizing if
\[
\int_{T}\int^{\sqrt{D(T)}}_0(\mu(B(t,r^2)))^{-{1}/{2}}\,dr \mu
(dt)<\infty.
\]
Let
\[
\ccM=\sup_{\mu} \int_{T}\int^{\sqrt{D(T)}}_0(\mu
(B(t,r^2)))^{-
{1}/{2}}\,dr \mu(dt),
\]
where the supremum is taken over all probability measures on $T$.

\begin{theo}[~\cite{Fer2,Tal1}]\label{thm1}
If $\ccM<\infty$, that is, all probability measures are weakly
majorizing with a uniform bounding constant, then there exists $m$ a
majorizing measure on $T$ such that
\[
\sup_{t\in T}\int^{\sqrt{D(T)}}_0 (m(B(t,r^2)))^{-{1}/{2}}\,dr
\ls\ccM.
\]
\end{theo}

The main result of this paper is the following theorem.

\begin{theo}\label{thm3}
Whenever all orthogonal processes on $T$ satisfy (\ref{square}), then
$\ccM\ls K D(T)$, where $K<\infty$.\vadjust{\goodbreak}
\end{theo}

When combined with Theorems~\ref{thm0},~\ref{thm1},~\ref{thm3} it implies
the complete characterization of a.e. convergence of all orthogonal series.

\begin{coro}
For a given $(a_n)^{\infty}_{n=1}$ series $\sum^{\infty
}_{n=1}a_n\varphi_n$
are a.e. convergent for all orthonormal sequences $(\varphi_n)^{\infty
}_{n=1}$ if and only if there exists a majorizing measure~$m$ on $T$.
\end{coro}

We stress that using the chaining argument and the Fernique's
idea of constructing a majorizing measure
makes the proof significantly shorter than the one presented in~\cite{Pas1}.

\section{Structure of the proof}\label{sect15}

If all orthogonal process satisfy (\ref{square}), then in particular
$D(T)<\infty$. For simplicity assume that $T\subset[0,1)$
(the general result can be obtained by the translation invariance and
homogeneity). Our approach is based on proving special properties
of natural partitions of $[0,1)$. Let
%
\begin{equation}\label{partia}
\qquad\ccA_{k}=\bigl\{A_i^{(k)}\dvtx  0\ls i<4^k\bigr\},\qquad k\gs0 \mbox{ where }
A^{(k)}_i=\bigl[i4^{-k},(i+1)4^{-k}\bigr)\cap T,
\end{equation}
in particular $A^{(0)}_0=T$. In Section~\ref{sect1} we translate the
weakly majorizing measure functionals into the language
of $\ccA_k$, $k\gs0$. Since as sated in Theorem~\ref{thm1} we have to
deal with any probability measure~$\mu$ on $T$, we fix~$\mu$ and check that for the particular $0\ls i<4^k$ sets
$A^{(k)}_{4i+j}$, $j\in\{0,1,2,3\}$, are important only if the measure~$\mu$ of $A^{(k-1)}_i$ is well distributed
among them. In this way we obtain the quantity that one may use to
bound the weakly majorizing measure functional.

Then we follow the idea that was first invented by Talagrand
in~\cite{Tal0} to prove
the complete characterization of Gaussian sample boundedness. We
introduce the set functionals
$F_k$, $k\gs0$, such that $F_k$ operates on $\ccA_k$ and is given by
%
\begin{equation}\label{functo}
F_k\bigl(A^{(k)}_{i}\bigr)=\sup_Y \E\sup_{t\in A^{(k)}_i} Y(t),
\end{equation}
where the supremum is over the class of processes $Y(t)$, $t\in\bar
{A}^{(k)}_i$, where $\bar{A}^{(k)}_i=A^{(k)}_i\cup\{
i4^{-k},(i+1)4^{-k}\}$,
that satisfy $\E Y(t)=0$ and
%
\begin{equation}\label{amilo}
\E|Y(s)-Y(t)|^2=|s-t|(1-4^{k}|s-t|)\qquad \mbox{for } s,t\in\bar{A}^{(k)}_i.
\end{equation}
In particular $Y(i4^{-k})=Y((i+1)4^{-k})$, and hence we may require
$Y(i4^{-k})=Y((i+1)4^{-k})=0$ [it does not change $F_k(A^{(k)}_i)$].
We show in Section~\ref{sect2} that if (\ref{square}) holds for all
orthogonal processes, then $F_0(T)<\infty$.
The partitioning scheme is the induction step which shows that
partitioning of $A^{(k-1)}_i$ into $A^{(k)}_{4i+j}$, $j\in\{0,1,2,3\}$,
makes it possible to earn the suitable quantity so that summing all
over the partitions completes the argument of
the uniform bound existence for any weakly majorizing measure
functional. The proof of the induction step\vspace*{1pt} is the construction
for a fixed $0\ls i< 4^{k-1}$ of a special process $Y(t)$, $t\in\bar
{A}^{k-1}_i$, that satisfies~(\ref{amilo}).
In the construction we use optimal (or nearly optimal) processes on
$A^{(k)}_{4i+j}$ for $j\in\{0,1,2,3\}$
and a suitably chosen family of independent random variables.

\section{Weakly majorizing measures}\label{sect1}

We have noted in Section~\ref{sect15} that one may assume $T\subset[0,1)$.
Consequently $\mu$ is weakly majorizing if
%
\begin{equation}\label{unt1}
\int_T \int^{1}_0 (\mu(B(t,r^2)))^{-{1}/{2}}\,dr \mu(dt)<\infty.
\end{equation}
We first translate the functional from (\ref{unt1}) into the language
of $\ccA_k$, $k\gs0$, defined in~(\ref{partia}).
\begin{lema}\label{lem1}
For each measure $\mu$ the inequality holds
\[
\int_T \int^{1}_0 (\mu(B(t,r^2)))^{-{1}/{2}}\,dr \mu(dt)\ls
\sum^{\infty}_{k=1}2^{-k}\sum^{4^k-1}_{i=0}\bigl(\mu\bigl(A^{(k)}_i\bigr)\bigr)^{{1}/{2}}.
\]
\end{lema}
\begin{pf}
First observe that
\[
\int^1_0 (\mu(B(t,r^2)))^{-{1}/{2}}\,dr\ls
\sum^{\infty}_{k=1}2^{-k}(\mu(B(t,4^{-k})))^{-{1}/{2}}\qquad \mbox{for } t\in T.
\]
Clearly $|\ccA_k|\ls4^{k}$ and $A^{(k)}_i\subset B(t,4^{-k})$ for all
$t\in A^{(k)}_i\in\ccA_k$. Consequently
$\mu(A^{(k)}_i)\ls\mu(B(t,4^{-k}))$, and hence
\begin{eqnarray*}
\int_T (\mu(B(t,4^{-k})))^{-{1}/{2}}\mu(dt)&\ls&
\sum^{4^k-1}_{i=0}\int_{A^{(k)}_{i}}(\mu(B(t,4^{-k})))^{-{1}/{2}}\mu
(dt)\\
&\ls&\sum^{4^k-1}_{i=0}\int_{A^{(k)}_i}\bigl (\mu\bigl(A^{(k)}_{i}\bigr)\bigr)^{-{1}/{2}}\mu(dt)
=\sum^{4^k-1}_{i=0}\bigl(\mu\bigl(A^{(k)}_i\bigr)\bigr)^{{1}/{2}}.
\end{eqnarray*}
Therefore
\begin{eqnarray*}
\int_T \int^{1}_0 (\mu(B(t,r)))^{-{1}/{2}}\,dr \mu(dt)&\ls&
\sum^{\infty}_{k=1}2^{-k}\int_T(\mu(B(t,4^{-k})))^{-
{1}/{2}}\mu
(dt)\\
&\ls&\sum^{\infty}_{k=1}2^{-k}\sum^{4^k-1}_{i=0}\bigl(\mu
\bigl(A^{(k)}_i\bigr)\bigr)^{{1}/{2}}.
\end{eqnarray*}
\upqed\end{pf}
For a specific measure $\mu$ not all subsets $A^{(k)}_i\in\ccA_k$ are
important.
Observe that for $0\ls i< 4^{k-1}$, $\bigcup
^3_{j=0}A^{(k)}_{4i+j}=A^{(k-1)}_{i}$.
Denote by $I(k)$ the set of indices $4i+j$ where $0\ls i<4^k$, $0\ls
j\ls3$ such that
%
\begin{equation}\label{mydlak01}
\tfrac{1}{32}\mu\bigl(A^{(k-1)}_{i}\bigr)\ls\mu\bigl(A^{(k)}_{4i+j}\bigr) \ls\tfrac
{1}{2}\mu
\bigl(A^{(k)}_{4i}\cup A^{(k)}_{4i+2}\bigr)
\end{equation}
if $j\in\{0,2\}$, and
%
\begin{equation}\label{mydlak02}
\tfrac{1}{32}\mu\bigl(A^{(k-1)}_{i}\bigr)\ls\mu\bigl(A^{(k)}_{4i+j}\bigr) \ls\tfrac
{1}{2}\mu
\bigl(A^{(k)}_{4i+1}\cup A^{(k)}_{4i+3}\bigr)
\end{equation}
if $j\in\{1,3\}$. The meaning of the construction is that $4i+j\in
I(k)$ only if measure
of $A^{(k-1)}_{i}$ is well distributed among $A^{(k)}_{4i+j}$, $j\in\{
0,1,2,3\}$.

We improve Lemma~\ref{lem1}, showing that the upper bound for the
weakly majorizing measure functional can be replaced by the one that
uses only sets of the form
$A^{(k)}_i$, $i\in I(k)$.
\begin{prop}\label{thm8}
For each probability Borel measure $\mu$ on $T$, the following
inequality holds:
\[
\int_T \int^1_0 (\mu(B(t,\va)))^{-{1}/{2}}\ls\frac
{1}{1-2^{-1}L}\Biggl[L+\sum_{k=1}^{\infty}2^{-k}\sum^{4^k-1}_{i=0}\bigl(\mu
\bigl(A^{(k)}_{i}\bigr)\bigr)^{{1}/{2}}1_{i\in I(k)}\Biggr],
\]
where $L=2^{{1}/{2}}\cdot\frac{5}{4}<2$.
\end{prop}
\begin{pf}
Suppose that $4i+j\notin I(k)$ and $j\in\{0,2\}$, then there are two
possibilities, either
%
\begin{eqnarray}
\label{set1}  \mu\bigl(A^{(k)}_{4i+j}\bigr)&<& \tfrac{1}{32}\mu\bigl(A^{(k-1)}_{i}\bigr),
\quad\mbox{or} \\
\label{set2}  \mu\bigl(A^{(k)}_{4i}\cup A^{(k)}_{4i+2}\bigr)&>&2\mu\bigl(A^{(k)}_{4i+j}\bigr).
\end{eqnarray}
If (\ref{set1}) holds, then
%
\begin{equation}\label{mydlak1}
\bigl(\mu\bigl(A^{(k)}_{4i+j}\bigr)\bigr)^{{1}/{2}}< \tfrac{2^{{1}/{2}}}{8}\bigl(\mu
\bigl(A^{(k-1)}_{i}\bigr)\bigr)^{{1}/{2}}.
\end{equation}
Assuming (\ref{set2}) we use the trivial inequality
%
\begin{equation}\label{mydlak2}
\bigl(\mu\bigl(A^{(k)}_i\bigr)\bigr)^{{1}/{2}}< \bigl(\mu\bigl(A^{(k)}_{4i}\cup
A^{(k)}_{4i+2}\bigr)\bigr)^{{1}/{2}}.
\end{equation}
One cannot have that both $j=0$ and $j=2$ satisfy (\ref{set2}), and
therefore due to (\ref{mydlak1}) and (\ref{mydlak2}),
%
\begin{eqnarray}\label{set3}
&&\bigl (\mu\bigl(A^{(k)}_{4i}\bigr)\bigr)^{{1}/{2}}1_{4i\notin I(k)}+\bigl(\mu
\bigl(A^{(k)}_{4i+2}\bigr)\bigr)^{{1}/{2}}1_{4i+2\notin I(k)}
\nonumber
\\
&&\qquad\ls\max\bigl\{\tfrac{2^{{1}/{2}}}{4}\bigl(\mu
\bigl(A^{(k-1)}_{i}\bigr)\bigr)^{{1}/{2}},\tfrac{2^{{1}/{2}}}{8}\bigl(\mu
\bigl(A^{(k-1)}_{i}\bigr)\bigr)^{{1}/{2}}\\
&&\hspace*{99pt}\qquad{}+\bigl(\mu\bigl(A^{(k)}_{4i}\cup
A^{(k)}_{4i+2}\bigr)\bigr)^{{1}/{2}}\bigr\}.\nonumber
\end{eqnarray}
The same argument works for $j\in\{1,3\}$, and consequently
%
\begin{eqnarray}\label{set4}
&& \bigl(\mu\bigl(A^{(k)}_{4i+1}\bigr)\bigr)^{{1}/{2}}1_{4i+1\notin I(k))}+\bigl(\mu
\bigl(A^{(k)}_{4i+3}\bigr)\bigr)^{{1}/{2}}1_{4i+3\notin I(k)}
\nonumber
\\
 &&\qquad \ls\max\bigl\{\tfrac{2^{{1}/{2}}}{4}\bigl(\mu
\bigl(A^{(k-1)}_{i}\bigr)\bigr)^{{1}/{2}},\tfrac{2^{{1}/{2}}}{8}\bigl(\mu
\bigl(A^{(k-1)}_{i}\bigr)\bigr)\\
&&\hspace*{78pt}\qquad{}+\bigl(\mu\bigl(A^{(k)}_{4i+1}\cup A^{(k)}_{4i+3}\bigr)\bigr)^{
{1}/{2}}\bigr\}.\nonumber
\end{eqnarray}
Since $x^{{1}/{2}}+y^{{1}/{2}}\ls2^{{1}/{2}}(x+y)^{{1}/{2}}$, for $x,y\gs0$ we have
%
\begin{equation}\label{legis3}
\quad\bigl(\mu\bigl(A^{(k)}_{4i}\cup A^{(k)}_{4i+2}\bigr)\bigr)^{{1}/{2}}+\bigl(\mu
\bigl(A^{(k)}_{4i+1}\cup A^{(k)}_{4i+3}\bigr)\bigr)^{{1}/{2}}\ls2^{
{1}/{2}}\bigl(\mu\bigl(A^{(k-1)}_{i}\bigr)\bigr)^{{1}/{2}}.
\end{equation}
On the other hand,
%
\begin{equation}\label{legis4}
\qquad\max\bigl\{\bigl(\mu\bigl(A^{(k)}_{4i}\cup A^{(k)}_{4i+2}\bigr)\bigr)^{{1}/{2}},\bigl(\mu
\bigl(A^{(k)}_{4i+1}\cup A^{(k)}_{4i+3}\bigr)\bigr)^{{1}/{2}}\bigr\}\ls\bigl(\mu
\bigl(A^{(k-1)}_{i}\bigr)\bigr)^{{1}/{2}}.
\end{equation}
By (\ref{legis3}) and (\ref{legis4}) we obtain that
\[
\sum^{3}_{j=0}\bigl(\mu\bigl(A^{(k)}_{4i+j}\bigr)\bigr)^{{1}/{2}}I_{4i+j\notin
I(k)}\ls L\bigl(\mu\bigl(A^{(k-1)}_i\bigr)\bigr),
\]
where $L=2^{{1}/{2}}\cdot\frac{5}{4}$. Consequently,
%
\begin{equation}\label{mydlak3}
\sum^{4^k-1}_{i=0} \bigl(\mu\bigl(A^{(k)}_i\bigr)\bigr)^{{1}/{2}}1_{i\notin
I(k)}\ls
L \sum^{4^{k-1}-1}_{i=0}\bigl(\mu\bigl(A^{(k-1)}_{i}\bigr)\bigr).
\end{equation}
Using (\ref{mydlak3}), we deduce
\begin{eqnarray*}
&&\sum^{\infty}_{k=1}2^{-k}\sum^{4^k-1}_{i=0}\bigl(\mu\bigl(A^{(k)}_i\bigr)\bigr)^{
{1}/{2}}\\
&&\qquad\ls
\sum^{\infty}_{k=1}2^{-k}\sum^{4^k-1}_{i=0}\bigl(\mu\bigl(A^{(k)}_i\bigr)\bigr)^{
{1}/{2}}1_{i\in I(k)}+L \sum^{\infty}_{k=1}2^{-k}\sum
^{4^{k-1}-1}_{i=0}\bigl(\mu\bigl(A^{(k-1)}_i\bigr)\bigr)^{{1}/{2}}.
\end{eqnarray*}
Since $\mu(A^{(0)}_0)=1$, it implies that
\[
(1-2^{-1}L)\sum^{\infty}_{k=1}2^{-k}\sum^{4^k-1}_{i=0}\bigl(\mu
\bigl(A^{(k)}_i\bigr)\bigr)^{{1}/{2}}\ls
L+ \sum^{\infty}_{k=1}2^{-k}\sum^{4^k-1}_{i=0}\bigl(\mu
\bigl(A^{(k)}_i\bigr)\bigr)^{
{1}/{2}}1_{i\in I(k)}.
\]
To complete the proof it suffices to apply Lemma~\ref{lem1}.
\end{pf}

\section{The partitioning scheme}\label{sect2}

In this section we prove the main induction procedure. Recall that
$(F_k)_{k\gs0}$ are set functionals defined in (\ref{functo}).
We are going to show that
%
\begin{equation}\label{andy1}
\sup_X \Bigl(\E\sup_{t\in T}\bigl(X(t)-X(0)\bigr)^2\Bigr)^{1/2}\gs\frac{1}{64}\sum
^{\infty
}_{k=0}2^{-k}\sum^{3}_{j=0}\bigl(\mu\bigl(A^{(k)}_{4i+j}\bigr)\bigr)^{
{1}/{2}}1_{4i+j\in I(k)},\hspace*{-35pt}
\end{equation}
where the supremum is taken over all orthogonal processes on $T$.
The idea of the proof is to first show that $F_0(T)<\sup_{X}(\E\sup
_{t\in T}(X(t)-X(0))^2)^{1/2}$. Then we establish the induction step so
that $(\mu(A^{(k-1)}_i))^{1/2}F_{k-1}(A^{(k-1)}_{i})$ can be used to
bound $\sum^3_{j=0} (\mu(A^{(k)}_{4i+j}))^{1/2}F_k(A^{(k)}_{4i+j})$
for all $k\gs1$ and $0\ls i<4^{k-1}$
together with some additional term required to get (\ref{andy1}).

First consider the special case of $A^{(0)}_0=T$. For each $Y(t)$,
$t\in\bar{A}^{(0)}_0$ satisfying~(\ref{amilo})
for $k=0$, we take $Z$ independent of $Y$ such that $\E Z=0$, $\E
Z^2=1$. Then the process
\[
X(t)=Y(t)+t Z,\qquad t\in T,
\]
satisfies (\ref{ineq2}) and, moreover, by Jensen's inequality,
%
\begin{equation}\label{bbb}
\E\sup_{t\in T}Y(t)=\E\sup_{t\in t}\bigl(Y(t)-Y(0)\bigr)\ls\Bigl(\E\sup_{t\in T}
\bigl(X(t)-X(0)\bigr)^2\Bigr)^{{1}/{2}}.
\end{equation}
Therefore (\ref{square}) implies that $F_0(T)<\infty$, which makes the
induction accessible.

The crucial idea is to show that the induction step is valid.
\begin{prop}\label{thm5}
For each $A^{(k-1)}_i$, $0\ls i< 4^{k-1}$ and $k\gs1$, the following
inequality holds:
\begin{eqnarray*}
&& \bigl(\mu\bigl(A^{(k-1)}_i\bigr)\bigr)^{{1}/{2}}F_{k-1}\bigl(A^{(k-1)}_{i}\bigr)\\
 &&\qquad\gs
\frac{1}{64}2^{-k}\sum^{3}_{j=0}\bigl(\mu\bigl(A^{(k)}_{4i+j}\bigr)\bigr)^{
{1}/{2}}1_{4i+j\in I(k)}+\sum^{3}_{j=0}\bigl(\mu\bigl(A^{(k)}_{4i+j}\bigr)\bigr)^{{1}/{2}}F_{k}\bigl(A^{(k)}_{4i+j}\bigr).
\end{eqnarray*}
\end{prop}
\begin{pf}
Fix $A^{(k-1)}_i$, $0\ls i<4^{k-1}$, $k\gs1$. We may assume that $\mu
(A^{(k-1)}_i)>0$, since otherwise there is nothing to prove.
On each $\bar{A}^{(k)}_{4i+j}$, $0\ls j\ls3$, there exist a process
$Y_j$, such that
\[
\E|Y_l(t)-Y_l(s)|^2=|t-s|(1-4^{k}|t-s|)\qquad \mbox{for } s,t\in\bar
{A}^{(k)}_{4i+j}
\]
and
%
\begin{equation}\label{a2}
\E\sup_{t\in A^{(k)}_{4i+j}}Y_j(t)\gs F_{k}\bigl(A^{(k)}_{4i+j}\bigr)-\va.
\end{equation}
As we have mentioned, we may assume that
$Y_j((4i+j)4^{-k})=Y_j((4i+j+1)4^{-k})=0$.
Our goal is to construct a process, $Y(t)$, $t\in\bar
{A}^{(k-1)}_{i}$, using
$Y_j$, $0\ls j\ls3$, that verifies (\ref{amilo}) for $\bar{A}^{(k-1)}_i$.

To construct $Y(t)$, $t\in T$, we will need also a family of independent
random variables $Z_j$, $0\ls j\ls3$.
We require that $Z_j$ are independent of processes $Y_j$, $0\ls j\ls3$,
and such that $\E Z_j=0$ and $\E Z_j^2=1$.
Let $S_0=0$ and for $1\ls j\ls4$,
\[
S_j=\sum^{j-1}_{l=0}Z_l- j 4^{-1}\Biggl(\sum^{3}_{l=0}Z_l\Biggr)\qquad \mbox{for }
1\ls j\ls4.
\]
Observe that for $0\ls l,m \ls4$,
\[
\E|S_{l}-S_{m}|^2=|l-m|(1-4^{-1}|l-m|).
\]
With the family $Z_j$, $0\ls j\ls3$, we associate a random variable
$\tau$
valued in $\{0,1,2,3\}$. We require that $\tau$ is independent of
$Y_j$, $0\ls j\ls3$,
and distributed as follows:
%
\begin{equation}\label{distr}
\P(\tau=j)=\frac{\mu(A^{(k)}_{4i+j})}{\mu(A^{(k-1)}_{i})}\qquad \mbox
{for } 0\ls j\ls3.
\end{equation}
We define the process $Y(t)$, $t\in\bar{A}^{(k)}_{4i+j}$, by
%
\begin{eqnarray}\label{proces}
Y(t)&=&2^{-k}S_{j}+2^{k}\bigl(t-(4i+j)4^{-k}\bigr)(S_{j+1}-S_{j})
\nonumber
\\[-8pt]
\\[-8pt]
\nonumber
&&{}+\bigl(\P(\tau
=j)\bigr)^{-{1}/{2}}Y_j(t)1_{\tau=j},
\end{eqnarray}
and also set $Y(i4^{-(k-1)})=Y((i+1)4^{-(k-1)})=0$.
We have to show that $Y(t)$, $t\in\bar{A}^{(k-1)}_{i}$, is admissible
for $F_k(A^{(k-1)}_i)$, that is, we make thorough calculations for the variance
of $Y(s)-Y(t)$, where $s,t\in\bar{A}^{(k-1)}_i$.
\begin{lema}\label{lem5}
The process $Y(t)$, $t\in\bar{A}^{(k-1)}_i$, satisfies $\E Y(t)=0$,
$t\in\bar{A}^{(k-1)}_i$, and
%
\begin{equation}\label{an3}
\E|Y(s)-Y(t)|^2=|s-t|(1-4^{k-1}|s-t|)\qquad \mbox{for } s,t\in\bar{A}^{(k-1)}_i.
\end{equation}
\end{lema}
\begin{pf}
The first assertion is trivial; we show (\ref{an3}).
Assume that $s,t\in\bar{A}^{(k)}_{4i+j}$, and then by (\ref{distr}),
the independence of $Z_j$, $0\ls j\ls3$, and
independence between $Z_j$, $0\ls j\ls3$, $\tau$ and $Y_j$, $0\ls
j\ls
3$ [recall that $\E Z_j=0$ and $\E Y_j(t)=0$,
$t\in\bar{A}^{(k)}_{4i+j}$]
we obtain that
\begin{eqnarray*}
&&\E|Y(s)-Y(t)|^2\\
&&\qquad=4^{k}|s-t|^2\E(S_{j+1}-S_{j})^2+\P(\tau=j)\P(\tau=j)^{-1}|s-t|(1-4^{k}|s-t|)\\
&&\qquad=
4^{k}(1-4^{-1})|t-s|^2+|s-t|(1-4^{k}|s-t|)=|s-t|(1-4^{k-1}|s-t|).
\end{eqnarray*}
Now suppose that $s\in\bar{A}^{(k)}_{4i+l}$, $t\in\bar
{A}^{(k)}_{4i+m}$ and $l< m$.
The idea we follow is to rewrite $|Y(s)-Y(t)|^2$ in terms of $Z_j$,
$0\ls j\ls3$ and $\tau$. Using that $Z_j$, $0\ls j\ls3$ are
independent and $Z_j$, $0\ls j\ls3$, $\tau$
are independent of $Y_j$, $0\ls j\ls3$
[moreover $\E Z_j=0$ and $\E Y_j(t)=0$, $t\in\bar{A}_{4i+j}$]
%
\begin{eqnarray}\label{ab0}
\E\bigl(Y(s)-Y(t)\bigr)^2&=&\E(Y_l(s))^2+\E(Y_m(t))^2\nonumber\\
&&{} +\E\bigl(Y(s)-\bigl(\P(\tau=l)\bigr)^{-{1}/{2}}Y_l(s)1_{\tau
=l}-Y(t)\\
&&\hspace*{71pt}{}+\bigl(\P(\tau=l)\bigr)^{-{1}/{2}}Y_m(t)1_{\tau=m}\bigr)^2.\nonumber
\end{eqnarray}
Clearly,
%
\begin{eqnarray}\label{ab1}
\E(Y_l(s))^2&=&\E
\bigl(Y_l(s)-Y_l\bigl((4i+l+1)4^{-k}\bigr)\bigr)^2
\nonumber
\\[-8pt]
\\[-8pt]
\nonumber
&=&|s-(4i+l+1)4^{-k}|\bigl(1-4^{-k}|s-(4i+l+1)4^{-k}|\bigr)
\end{eqnarray}
and
%
\begin{eqnarray}\label{ab2}
\E(Y_m(t))^2&=&\E
\bigl(Y_m\bigl((4i+m)4^{-k}\bigr)-Y_j(t)\bigr)^2
\nonumber
\\[-8pt]
\\[-8pt]
\nonumber
&=&|(4i+m)4^{-k}-t|\bigl(1-4^{-k}|(4i+m)4^{-k}-t|\bigr).
\end{eqnarray}
Then we observe that by the definition,
\begin{eqnarray*}
Y(s)-\bigl(\P(\tau=l)\bigr)^{-{1}/{2}}Y_l(s)1_{\tau
=l}&=&2^{-k}S_{l}+2^{k}\bigl(s-(4i+l)4^{-k}\bigr)(S_{l+1}-S_{l}),\\
Y(t)-\bigl(\P(\tau=m)\bigr)^{-{1}/{2}}Y_l(s)1_{\tau
=m}&=&2^{-k}S_{m}+2^{k}\bigl(t-(4i+m)4^{-k}\bigr)(S_{m+1}-S_{m}).
\end{eqnarray*}
Hence
\begin{eqnarray*}
&&Y(s)-\bigl(\P(\tau=l)\bigr)^{-{1}/{2}}Y_l(s)1_{\tau=l}-Y(t)+\bigl(\P(\tau
=l)\bigr)^{-{1}/{2}}Y_m(t)1_{\tau=m}\\
&&\qquad=2^{-k}(S_{m}-S_{l})\\
&&\qquad\quad{}+2^{k}\bigl[\bigl(t-(4i+m)4^{-k}\bigr)(S_{m+1}-S_{m})-\bigl(\bigl(s-(4i+l)4^{-k}\bigr)\bigr)(S_{l+1}-S_l)\bigr].
\end{eqnarray*}
Since $S_j=\sum^{j-1}_{l=0}Z_l- j 4^{-1}(\sum^{3}_{l=0}Z_l)$, we have
\begin{eqnarray*}
&&
2^{-k}(S_{m}-S_{l})\\
&&\quad{}+2^{k}\bigl[\bigl(t-(4i+m)4^{-k}\bigr)(S_{m+1}-S_{m})-\bigl(\bigl(s-(4i+l)4^{-k}\bigr)\bigr)(S_{l+1}-S_l)\bigr]\\
&&\qquad=2^{-k}\Biggl(\sum^{m-1}_{j=l}Z_j-\frac{m-j}{4}\sum
^3_{j=0}Z_j\Biggr)\\
&&\qquad\quad{}+2^{k}\Biggl[\bigl(t-(4i+m)4^{-k}\bigr)\Biggl(Z_m-\frac{1}{4}\sum
^3_{j=0}Z_j\Biggr)\Biggr]\\
&&\quad\qquad{}-2^{k}\Biggl[\bigl(\bigl(s-(4i+l)4^{-k}\bigr)\bigr)\Biggl(Z_l-\frac{1}{4}\sum^3_{j=0}Z_j\Biggr)\Biggr].
\end{eqnarray*}
We group coefficients by random variables $Z_j$. For $Z_m$ we obtain
\begin{eqnarray*}
&&-2^{-k}\frac{m-l}{4}+2^k\bigl(t-(4i+m)4^{-k}\bigr)\frac
{3}{4}+2^{k}\bigl(s-(4i+l)4^{-k}\bigr)\frac{1}{4}\\
&&\qquad=2^k|(4i+m)4^{-k}-t|-4^{-1}\bigl(2^{-k}(m-l)+2^k|s-t|-2^{-k}(m-l)\bigr)\\
&&\qquad=2^k\bigl(|(4i+m)4^{-k}-t|-4^{-1}|s-t|\bigr).
\end{eqnarray*}
Similarly the coefficient for $Z_l$ equals
\begin{eqnarray*}
&&-2^{-k}\frac{m-l}{4}-2^k\bigl(s-(4i+l)4^{-k}\bigr)\frac
{3}{4}-2^{k}\bigl(t-(4i+m)4^{-k}\bigr)\frac{1}{4}\\
&&\qquad=2^k|(4i+l+1)4^{-k}-s|\\
&&\qquad\quad{}-4^{-1}\bigl(2^{-k}(m-l-1)+2^k|s-t|-2^{-k}(m-l-1)\bigr)\\
&&\qquad=2^k\bigl(|s-(4i+l+1)4^{-k}|-4^{-1}|s-t|\bigr)2^k.
\end{eqnarray*}
For $l<j<m$ the coefficient for $Z_j$ is
\begin{eqnarray*}
&&2^{-k}\biggl(1-\frac{m-l}{4}\biggr)-2^k\bigl(t-(4i+m)4^{-k}\bigr)\frac
{1}{4}+2^k\bigl(s-(4i+l)4^{-k}\bigr)\frac{1}{4}\\
&&\qquad=2^{k}\bigl(4^{-k}-4^{-1}(m-l)4^{-k}-4^{-1}\bigl(|s-t|-(m-l)4^{-k}\bigr)\bigr)\\
&&\qquad=2^k(4^{-k}-4^{-1}|s-t|)
\end{eqnarray*}
and finally for $j>m$ and $j<l$
\[
-2^{-k}\frac{m-l}{4}-2^k\bigl(t-(4i+m)4^{-k}\bigr)\frac
{1}{4}+2^k\bigl(s-(4i+l)4^{-k}\bigr)\frac{1}{4}=-2^k(4^{-1}|s-t|).
\]
Consequently we obtain that
\begin{eqnarray*}
&&Y(s)-\bigl(\P(\tau=l)\bigr)^{-{1}/{2}}Y_l(s)1_{\tau=l}-Y(t)+\bigl(\P(\tau
=l)\bigr)^{-{1}/{2}}Y_m(t)1_{\tau=m}\\
&&\qquad=\bigl(|(4i+m)4^{-k}-t|-4^{-1}|s-t|\bigr)2^k
Z_m\\
&&\qquad\quad{}+\bigl(|s-(4i+l+1)4^{-k}|-4^{-1}|s-t|\bigr)2^k Z_l\\
&&\qquad\quad{}+\sum^{m-1}_{n=l+1}(4^{-k}-4^{-1}|s-t|)2^k Z_n-4^{-1}|s-t|2^k\sum
_{n<l,n>m}Z_n.
\end{eqnarray*}
Therefore by the orthogonality of $Z_j$, $j\in\{0,1,2,3\}$,
%
\begin{eqnarray}\label{ab3}
&& \E\bigl(Y(s)-\bigl(\P(\tau=l)\bigr)^{-{1}/{2}}Y_l(s)1_{\tau=l}-Y(t)+\bigl(\P
(\tau
=l)\bigr)^{-{1}/{2}}Y_m(t)1_{\tau=m}\bigr)^2\nonumber\\
&&\qquad=\bigl(|(4i+m)4^{-k}-t|-4^{-1}|s-t|\bigr)^2
4^k
\nonumber
\\[-8pt]
\\[-8pt]
\nonumber
&&\qquad\quad{}+\bigl(|s-(4i+l+1)4^{-k}|-4^{-1}|s-t|\bigr)^2 4^k\\
&&\qquad\quad{}+(4^{-k}-4^{-1}|s-t|)^2(m-l-1) 4^k+ 4^{-2}|s-t|^2
(4-m+l-1) 4^{k}.\nonumber
\end{eqnarray}
Combining (\ref{ab0}), (\ref{ab1}), (\ref{ab2}), (\ref{ab3}) and
\[
|s-(4i+l+1)4^{-k}|+(m-l-1)4^{-k}+|(4i+m)4^{-k}-t|=|s-t|,
\]
we obtain that
\[
\E|Y(s)-Y(t)|^2=|s-t|(1-4^{k-1}|s-t|).
\]
This completes the proof.
\end{pf}

Having the process $Y(t)$, $t\in\bar{A}^{(k-1)}_i$, constructed, we
use it
to provide a lower bound on $F_{k-1}(A^{(k-1)}_i)$.
First note that
\[
F_k\bigl(A^{(k-1)}_i\bigr)\gs\E\sup_{t\in A^{(k-1)}_i}Y(t)\gs\sum^3_{j=0}\E
\Bigl(\sup
_{t\in A^{(k)}_{4i+j}} Y(t)1_{\tau=j}\Bigr).
\]
Moreover,
\begin{eqnarray*}
&&\E\Bigl(\sup_{t\in A^{(k)}_{4i+j}} Y(t)1_{\tau=j}\Bigr)\\[-2pt]
&&\qquad =2^{-k}\E S_{j}1_{\tau =j}+\E\Bigl(\sup_{t\in
A^{(k)}_{4i+j}}\bigl(2^{k}\bigl(t-(4i+j)4^{-k}\bigr)(S_{j+1}-S_{j})1_{\tau=j}\\[-2pt]
&&\hspace*{209pt}{}+\bigl(\P
(\tau
=j)\bigr)^{-{1}/{2}}Y_j(t)\bigr)1_{\tau=j}\Bigr).
\end{eqnarray*}
Conditioning on $\ccF=\sigma(Y_j, 0\ls j\ls3)$ and then using Jensen's
inequality, we deduce
\begin{eqnarray*}
&&\E\sup_{t\in A^{(k)}_{4i+j}}\bigl(\bigl(2^{k}\bigl(t-(4i+j)\bigr)(S_{j+1}-S_{j})1_{\tau
=j}+\bigl(\P(\tau=j)\bigr)^{-{1}/{2}}Y_j(t)\bigr)1_{\tau=j}\bigr)\\[-2pt]
&&\qquad\gs\E\sup_{t\in A^{(k)}_{4i+j}}\bigl(\E
\bigl(\bigl(2^{k}\bigl(t-(4i+j)4^{-k}\bigr)(S_{j+1}-S_{j})1_{\tau=j}\\[-2pt]
&&\hspace*{127pt}\qquad{}+\bigl(\P(\tau
=j)\bigr)^{-{1}/{2}}Y_j(t)\bigr)1_{\tau=j}|\ccF\bigr)\bigr)\\[-2pt]
&&\qquad=\E\sup_{t\in A^{(k)}_{4i+j}}\bigl(2^{k}\E
\bigl(\bigl(t-(4i+j)4^{-k}\bigr)(S_{j+1}-S_{j})1_{\tau=j}\bigr)+
\bigl(\P(\tau=j)\bigr)^{{1}/{2}}Y_j(t)\bigr)\\[-2pt]
&&\qquad\gs-2^{-k}\bigl(\E(S_{j+1}-S_{j})1_{\tau=j}\bigr)_{-}+\bigl(\P(\tau=j)\bigr)^{{1}/{2}}
\E\sup_{t\in A^{(k)}_{4i+j}}Y_j(t).
\end{eqnarray*}
Consequently,
\begin{eqnarray*}
F_{k-1}\bigl(A^{(k-1)}_{i}\bigr)&\gs&\sum^3_{j=0}\Bigl(2^{-k}\bigl[\E S_{j}1_{\tau=j}-\bigl(\E
(S_{j+1}-S_{j})1_{\tau=j}\bigr)_{-}\bigr]\\[-2pt]
&&\hspace*{38pt}\qquad{}+\bigl(\P(\tau=j)\bigr)^{{1}/{2}}
\E\sup_{t\in A^{(k)}_{4i+j}}Y_j(t)\Bigr).
\end{eqnarray*}
Together with (\ref{a2}) and (\ref{distr}) it implies that
%
\begin{eqnarray}\label{bumbum}
F_{k-1}\bigl(A^{(k-1)}_i\bigr)&\gs&\sum^3_{j=0}\bigl(2^{-k}\bigl[\E S_{j}1_{\tau=j}-\bigl(\E
(S_{j+1}-S_{j})1_{\tau=j} \bigr)_{-}\bigr]
\nonumber
\\[-8pt]
\\[-8pt]
\nonumber
&&\hspace*{18pt}{}+\bigl(\P(\tau=j)\bigr)^{
{1}/{2}}\bigl(F_{k}(A^{k}_{4i+j})-4\va\bigr).
\end{eqnarray}
To complete the lower bound, we have to construct variables $Z_j$,
$0\ls j\ls3$, and~$\tau$. The main idea is to choose $n\in\{0,1,2,3\}$
and variable $Z_n$
to be $\tau$ measurable, whereas all remaining $Z_j$, $j\neq n$, are
independent of $\tau$. Therefore we first define $\tau$ so that (\ref
{distr}) holds, then
obtain $Z_n$ as a Borel function of $\tau$ and only then set any
independent $Z_j$, $j\neq n$,
independent of $Z_n$. In this setting, define
\[
V_n=\sum^3_{j=0}\bigl(\E S_{j}1_{\tau=j}-\bigl(\E(S_{j+1}-S_{j})1_{\tau=j} \bigr)_{-}\bigr).
\]
Observe that since $Z_l$, $l\neq n$, are independent of $\tau$ and
consequently of $Z_{n}$, we have
$\E Z_l 1_{\tau=j}=\E Z_l \P(\tau=j)=0$, whenever $l\neq n$. Therefore
\[
\E S_j 1_{\tau=j}=\E Z_n 1_{\tau=j} 1_{n\ls j-1}-\frac{j}{4}\E Z_n
1_{\tau=j}
\]
and
\[
\E(S_{j+1}-S_j)1_{\tau=j}=\E Z_n 1_{\tau=j}1_{j=n}-\tfrac{1}{4}\E
Z_{n}1_{\tau=j}.
\]
Consequently for $j\neq n$, $(\E(S_{j+1}-S_j)1_{\tau=j})_{-}=-\frac
{1}{4}(\E Z_n 1_{\tau=j})_{+}$
and for $j=n$, $(\E(S_{j+1}-S_j)1_{\tau=j})_{-}=\frac{3}{4}(\E Z_n
1_{\tau=j})_{-}$.
Hence the representation
\[
V_n=\sum^3_{j=n+1}c_j-(1-4^{-1})(c_n)_{-}-\sum
^{3}_{j=0}j4^{-1}c_j-\sum
_{l\neq n}4^{-1}(c_j)_{+},
\]
where $c_j=\E Z_{n} 1_{\tau=j}$. Since $\va>0$ is arbitrary in (\ref
{bumbum}),
we obtain
%
\begin{equation}\label{bambam1}
F_{k-1}\bigl(A^{(k-1)}_i\bigr)\gs2^{-k}V_n+\sum^3_{j=0} \bigl(\P(\tau=j)\bigr)^{
{1}/{2}}F_{k}\bigl(A^{(k)}_{4i+j}\bigr).
\end{equation}
The above inequality completes the first part of the proof. Using the process
$Y(t)$, $t\in\bar{A}^{(k-1)}_i$, we have shown that $\mu
(A^{(k-1)}_i)F_{k-1}(A^{(k-1)}_i)$ dominates $\sum^3_{j=0}\mu
(A^{(k)}_{4i+j})F_k(A^{(k)}_{4i+j})$,\vspace*{2pt} together with the additional term
$2^{-k}V_n$.

We claim that it is always possible to define $Z_n$ with respect to
$\tau$
in a way that one can bound $V_n$ from below by a universal constant,
assuming that there exists at least one $j\in\{0,1,2,3\}$ such that
$4i+j\in I(k)$.
\begin{lema}\label{pip} There exists $Z_3$ measurable with respect to
$\tau$, such that $\E Z_3=0$, $\E Z_3^2=1$ and
%
\begin{equation}\label{muba1}
V_3\gs\frac{1}{4}\biggl(\frac{\P(\tau=0)\P(\tau=2)}{\P(\tau=0)+\P
(\tau
=2)}\biggr)^{{1}/{2}}
\end{equation}
and $Z_2$ measurable with respect to $\tau$, such that $\E Z_2=0$, $\E
Z_2^2=1$ and
%
\begin{equation}\label{muba2}
V_2\gs\frac{1}{4}\biggl(\frac{\P(\tau=1)\P(\tau=3)}{\P(\tau=1)+\P
(\tau
=3)}\biggr)^{{1}/{2}}.
\end{equation}
\end{lema}
\begin{pf}
First note\vspace*{1pt} that $\sum^3_{j=0}c_j=0$, and then observe that it benefits
to set $c_n=0$.
The first case we consider is $n=3$, so $c_3=0$, and
then if $c_0\gs0$, $c_1=0$, $c_2\ls0$, we have
%
\begin{equation}\label{ados1}
V_3=-\tfrac{1}{4}c_0-\tfrac{2}{4}c_2=-\tfrac{1}{4}c_2=\tfrac{1}{4}c_0,
\end{equation}
where we have used that $c_0+c_2=0$. The second case is when $n=2$,
$c_2=0$, and then if
$c_0=0$, $c_1\ls0$, $c_3\gs0$, we have
%
\begin{equation}\label{ados2}
V_2=c_3-\tfrac{1}{4}c_1-\tfrac{3}{4}c_3-\tfrac{1}{4}c_3=-\tfrac
{1}{4}c_1=\tfrac{1}{4}c_3,
\end{equation}
where we have used that $c_1+c_3=0$. In the same way one can treat
$V_0$ and~$V_1$.

The above discussion leads to the definition of $Z_n$. If
$n=3$, we set
\[
Z_3=x1_{\tau=0}+y1_{\tau=2}.
\]
Our requirements are $\E Z_3=0$, $\E Z_3^2=1$, so
\begin{eqnarray*}
x\P(\tau=0)+y\P(\tau=2)&=&0, \\
x^2\P(\tau=0)+y^2\P(\tau=2)&=&1.
\end{eqnarray*}
Therefore
\[
x=\biggl(\frac{\P(\tau=2)}{\P(\tau=0)(\P(\tau=0)+\P(\tau=2))}
\biggr)^{{1}/{2}},
\]
and consequently all the requirements for (\ref{ados1}) are satisfied,
and we have
\[
V_3=\frac{1}{4}c_0=\frac{1}{4}\biggl(\frac{\P(\tau=0)\P(\tau=2)}{\P
(\tau
=0)+\P(\tau=2)}\biggr)^{{1}/{2}}.
\]
The same argument for $n=2$ shows that one can construct $Z_2$ in a way that
all requirements for (\ref{ados2}) are satisfied and
\[
V_2=\frac{1}{4}c_3=\frac{1}{4}\biggl(\frac{\P(\tau=1)\P(\tau=3)}{\P
(\tau
=1)+\P(\tau=3)}\biggr)^{{1}/{2}}.
\]
\upqed\end{pf}
We use the above lemma in (\ref{bambam1}) to bound $2^{-k}V_n$.
There are three cases. First suppose that $4i+j\notin I(k)$ for $0\ls
j\ls3$, and then we set $Z_j$, $j\in\{0,1,2,3\}$, to be independent of
$\tau$ which implies that
$V_n= 0$ for any choice of $n$. Therefore by~(\ref{bambam1}),
\begin{eqnarray*}
F_{k-1}\bigl(A^{(k-1)}_i\bigr)&\gs&\sum^3_{j=0}\bigl(\P(\tau=j)\bigr)^{
{1}/{2}}F_k\bigl(A^{(k)}_{4i+j}\bigr)\\
&=&\frac{1}{64}2^{-k}\sum^{4i+3}_{j=4i}\bigl(\P(\tau=j)\bigr)^{
{1}/{2}}1_{4i+j\in I(k)}\\
&&{}+
\sum^3_{j=0}\bigl(\P(\tau=j)\bigr)^{{1}/{2}}F_k\bigl(A^{(k)}_{4i+j}\bigr).
\end{eqnarray*}
The second case is that $4i+j\in I(k)$ for $j\in\{0,2\}$, then we use
(\ref{mydlak01}) and (\ref{muba1})
\begin{eqnarray*}
V_3&\gs&\frac{1}{4}\biggl(\frac{\P(\tau=0)\P(\tau=2)}{\P(\tau=0)+\P
(\tau
=2)}\biggr)^{{1}/{2}}\gs\frac{1}{4}\biggl(\frac{1}{2}\cdot\frac
{1}{32}\biggr)^{
{1}/{2}}=\frac{1}{32} \\
&\gs&\frac{1}{64}\sum^{4i+3}_{j=4i}\bigl(\P(\tau
=j)\bigr)^{{1}/{2}}1_{j\in I(k)},
\end{eqnarray*}
where we have used the inequality $x^{{1}/{2}}+y^{
{1}/{2}}+z^{{1}/{2}}+t^{{1}/{2}}\ls2(x+y+z+t)^{
{1}/{2}}$, for
$x,y,z,t\gs0$. Therefore
\begin{eqnarray*}
F_{k-1}\bigl(A^{(k-1)}_i\bigr)&\gs&2^{-k}V_3+\sum^3_{j=0}\bigl(\P(\tau=j)\bigr)^{
{1}/{2}}F_k\bigl(A^{(k)}_{4i+j}\bigr)\\
&\gs&\frac{1}{64}2^{-k}\sum^{4i+3}_{j=4i}\bigl(\P(\tau=j)\bigr)^{
{1}/{2}}1_{4i+j\in I(k)}\\
&&{}+
\sum^3_{j=0}\bigl(\P(\tau=j)\bigr)^{{1}/{2}}F_k\bigl(A^{(k)}_{4i+j}\bigr).
\end{eqnarray*}
The third possibility is that $4i+j\in I(k)$, $j\in\{1,3\}$, and then
by (\ref{mydlak02}) and (\ref{muba2})
we have
\begin{eqnarray*}
V_2&\gs&\frac
{1}{4}\biggl(\frac{\P(\tau=1)\P(\tau=3)}{\P(\tau=1)+\P(\tau
=3)}\biggr)^{
{1}/{2}}\gs\frac{1}{4}\biggl(\frac{1}{2}\cdot\frac{1}{32}\biggr)^{{1}/{2}}=
\frac{1}{32}\\
&\gs&\frac{1}{64}\sum^{4i+3}_{j=4i}\bigl(\P(\tau=j)\bigr)^{
{1}/{2}}1_{4i+j\in I(k)}.
\end{eqnarray*}
Consequently
\begin{eqnarray*}
F_{k-1}\bigl(A^{(k-1)}_i\bigr)&\gs&2^{-k}V_3+\sum^3_{j=0}\bigl(\P(\tau=j)\bigr)^{
{1}/{2}}F_k\bigl(A^{(k)}_{4i+j}\bigr)\\
&\gs&\frac{1}{64}2^{-k}\sum^{4i+3}_{j=4i}\bigl(\P(\tau=j)\bigr)^{
{1}/{2}}1_{4i+j\in I(k)}\\
&&{}+
\sum^3_{j=0}\bigl(\P(\tau=j)\bigr)^{{1}/{2}}F_k\bigl(A^{(k)}_{4i+j}\bigr).
\end{eqnarray*}
In the view of (\ref{distr}) it completes the proof of Proposition
\ref{thm5}.
\end{pf}

\section{Proof of the main result}\label{sect3}

In this section we use the functional $F_k$, $k\gs0$, and the
induction scheme proved in Proposition~\ref{thm5} to prove (\ref{andy1}).
\begin{prop}\label{prot}
The following inequality holds:
\[
\sum^{\infty}_{k=1}2^{-k}\sum^{4^k-1}_{i=0}\bigl(\mu\bigl(A^{(k)}_i\bigr)\bigr)^{
{1}/{2}}1_{i\in I(k)}\ls64\Bigl(\sup_X\E\sup_{t\in T}
\bigl(X(t)-X(0)\bigr)^2\Bigr)^{{1}/{2}},
\]
where the supremum is taken over all orthogonal process on $T$.
\end{prop}
\begin{pf}
By (\ref{bbb}) we have
\[
F_0(T)\ls\Bigl(\sup_{X}\E\sup_{t\in T} \bigl(X(t)-X(0)\bigr)^2\Bigr)^{{1}/{2}}.
\]
On the other hand using the induction step proved in Proposition~\ref{thm5},
we deduce
\[
\sum^{\infty}_{k=1}2^{-k}\sum^{4^k-1}_{i=0}\bigl(\mu\bigl(A^{(k)}_i\bigr)\bigr)^{{1}/{2}}1_{i\in I(k)}\ls64 F_0(T).
\]
This completes the proof.
\end{pf}

Using Propositions~\ref{thm8} and~\ref{thm5}, we conclude Theorem
\ref
{thm3} with
\[
K=\frac{1}{1-2^{-1}L}\Bigl(L+64 \sup_X\Bigl(\E\sup_{t\in
T}\bigl(X(t)-X(0)\bigr)^2\Bigr)^{{1}/{2}}\Bigr),
\]
and $L=2^{{1}/{2}}\cdot\frac{5}{4}$.



\printaddresses

\end{document}